\documentclass[10pt]{article}

\usepackage{amsmath}
\usepackage{amsthm}
\usepackage{amssymb}
\usepackage{fullpage}
\usepackage{graphicx}
\usepackage{algorithmic}
\usepackage{algorithm}
\usepackage{verbatim}

\def\q{\mathbf{q}}
\def\Q{\mathbf{Q}}

\def\etal{et al.}
\def\ie{i.e.\ }

\begin{document}


\title{Asynchronous Variational Integration of Interaction Potentials
  for Contact Mechanics}

\author{Etienne Vouga \\{\footnotesize Columbia University} %
\and David Harmon \\{\footnotesize Columbia University} %
\and Rasmus Tamstorf\\{\footnotesize Walt Disney Animation Studios} %
\and Eitan Grinspun\\{\footnotesize Columbia University}}
\date{}
\maketitle

\section{Introduction}
\emph{Variational integrators} (VIs)~\cite{aSuris1990, cMacKay1992,
  aMarsden2001} are a general class of time integration methods for
Hamiltonian systems whose construction guarantees certain highly
desirable properties. Instead of directly discretizing the smooth
equations of motion of a system, the variational approach asks that we
instead step back and discretize the system's Lagrangian. By analogy
to Hamilton's Least Action Principle, we may then form a discrete
action and seek paths which extremize it, yielding \emph{discrete}
Euler-Lagrange equations from which discrete equations of motion are
readily recovered. As a consequence of this special, more principled
construction, variational integrators are guaranteed to satisfy a
discrete formulation of Noether's Theorem~\cite{pWest2004}, and as a
special case conserve linear and angular momentum. VIs are
automatically \emph{symplectic}~\cite{bHairer2006}; while they do not
necessarily conserve energy, conservation of the symplectic form
assures no-drift conservation of energy over exponentially many time
steps~\cite{bHairer2006}.

Mechanical systems are almost never uniformly stiff. Different
potentials have different stable time step requirements, and even for
identical potentials this requirement depends on element size, since
finer elements can support higher-energy modes than coarser
elements. Any global time integration scheme cannot take advantage of
this variability, and instead must integrate the entire system at the
globally stiffest time step. Suppose the system can be triangulated
into elements such that each force acts entirely within one
element. Then \emph{asynchronous variational integrators}
\cite{aLew2003} generalize VIs by allowing each element to have its
own, independent time step. Coarser elements can then be assigned a
slower ``clock,'' and finer elements a faster one, so that relatively
few very fine elements do not as significantly degrade the overall
performance of integrating the system. AVIs retain all of the
properties of variational integrators mentioned above, except for
symplecticity. However, AVIs instead preserve an analogous
\emph{multisymplectic} form, and it has been shown
experimentally that preservation of this form likely induces the same
long-time good energy behavior that characterize symplectic
integrators~\cite{aLew2003}.

The published proof of multisymplecticity assumes
that the potentials are of an ``elastic type,'' i.e., specified by
volume integration over the material domain, an assumption violated by
interaction-type potentials. We extend the proof, showing that AVIs
remain multisymplectic under relaxed assumptions on the type of
potential (\S\ref{sec:proof}). The modified proof allows for
interaction potentials of the kind needed for contact mechanics (i.e.,
penalty forces). The extended theory thus enables the simulation of
mechanical contact in elastica (such as thin shells) and multibody
systems (such as granular materials) with no drift of conserved
quantities (energy, momentum) over long run times, using the
algorithms in~\cite{aHarmon2009}.

We conclude with data from numerical experiments measuring the long
time energy behavior of simulated contact, comparing the method built
on multisymplectic integration of interaction potentials to recently
proposed methods for thin shell contact (\S\ref{sec:results}).

\section{Variational Integrators}\label{sec:proof}
We begin with a background on variational integration and symplectic
structure~\cite{bHairer2006, aMarsden2001, pWest2004}.

Let $\gamma(t)$ be a piecewise-regular trajectory through
configuration space $\Q$, and $\dot\gamma(t) = \frac{d}{dt}\gamma(t)$
be the configurational velocity at time $t$. For simplicity we shall
assume that the kinetic energy of the system $T$ depends only on
configurational velocity, and that the potential energy $V$ depends
only on configurational position, so that we may write the Lagrangian
$L$ at time $t$ as
\begin{align}
L(q, \dot q) = T(\dot q) - V(q). \label{eq-action}
\end{align}

Then given the configuration of the system $q_0$ at time $t_0$ and $q_f$ at $t_f$, Hamilton's principle~\cite{bLanczos1986} states that the trajectory of the system $\gamma(t)$ joining $\gamma(t_0) = q_0$ and $\gamma(t_f) = q_f$ is a stationary point of the action functional
\begin{align*}
S(\gamma) = \int_{t_0}^{t_f} L\left[ \gamma(t), \dot\gamma(t) \right] dt
\end{align*}
with respect to taking variations $\delta \gamma$ of $\gamma$ which leave $\gamma$ fixed at the endpoints $t_0,\, t_f$. In other words, $\gamma$ satisfies
\begin{align}
dS(\gamma) \cdot \delta \gamma = 0. \label{eq-Hamilton}
\end{align}

Integrating by parts, and using that $\delta \gamma$ vanishes at $t_0$ and $t_1$, we compute
\begin{align*}
dS(\gamma) \cdot \delta \gamma = \int_{t_0}^{t_f} \left( \frac{\partial L}{\partial q}(\gamma,\dot\gamma) \cdot \delta \gamma + \frac{\partial L}{\partial \dot q} (\gamma, \dot \gamma) \cdot \delta \dot\gamma \right) dt = \int_{t_0}^{t_f} \left( -\frac{\partial V}{\partial q}(\gamma) - \frac{\partial^2 T}{\partial \dot q^2} (\dot \gamma)\ddot\gamma \right) \cdot \delta \gamma \, dt = 0.
\end{align*}

Since this equality must hold for all variations $\delta \gamma$ that fix $\gamma$'s endpoints, we must have
\begin{align}
\frac{\partial V}{\partial q}(\gamma) + \frac{\partial^2 T}{\partial \dot q^2}(\dot \gamma)\ddot\gamma = 0 \label{eq-EL},
\end{align}
the \emph{Euler-Lagrange equation} of the system. This equation is a second-order ordinary differential equation, and so has a unique solution $\gamma$ given two initial values $\gamma(t_0)$ and $\dot \gamma(t_0)$.

\subsection{Symplecticity}

The flow $\Theta_s: \left[\gamma(t), \dot \gamma(t) \right] \mapsto \left[\gamma(t+s), \dot \gamma(t+s)\right]$ given by (\ref{eq-EL}) has many structure-preserving properties; in particular it is momentum-preserving, energy-preserving, and symplectic~\cite{tLew2003}. To see this last property, for the remainder of this section we restrict the space of trajectories to those that satisfy the Euler-Lagrange equations. For such trajectories, and relaxing the requirement that $\delta\gamma$ fix the endpoints of $\gamma$, we have
\begin{align}
dS(\gamma) \cdot \delta \gamma = \frac{\partial T}{\partial \dot q}\left[ \pi_{\dot q} (q, \dot q) \right]\cdot \delta \gamma\bigg|_{t_0}^{t_f},\label{eq-ELparam}
\end{align}
where $\pi_{\dot q}$ is projection onto the second factor.

Since initial conditions $(q, \dot q)$ are in bijection with
trajectories satisfying the Euler-Lagrange equation, such trajectories
$\gamma$ can be uniquely parameterized by initial conditions
$\left[\gamma(t_0), \dot\gamma(t_0)\right]$. For the remainder of this
section we also restrict variations $\delta \gamma$ to \emph{first variations}: those variations in whose direction
 $\gamma$ continues to satisfy the Euler-Lagrange equations. These
are also parameterized by variations of the initial conditions,
$(\delta q, \delta \dot q)$. For conciseness of notation, we will
write $\nu(t) = (\gamma(t), \dot \gamma(t))$ and $\delta \nu(t) =
\left[\delta \gamma(t), \delta \dot \gamma(t)\right]$; using this
notation we write the above two facts as $\nu(t) =
\Theta_{t-t_0}\nu(t_0)$ and $\delta\nu(t) = {\Theta_{t-t_0}}_*
\delta\nu(t_0)$. The action (\ref{eq-action}), a functional on
trajectories $\gamma$, can also be rewritten as a function $S_i$ of
the initial conditions,
\begin{align*}
S_i(q, \dot q) &= \int_0^{t_f-t_0} L\left[\Theta_t(q, \dot q) \right]dt,
\intertext{so that}
dS(\gamma) \cdot \delta \gamma &= dS_i\left[\nu(t_0) \right] \cdot \delta \nu(t_0).
\intertext{Substituting all of these expressions into (\ref{eq-ELparam}), we get}
dS_i\left[\nu(t_0)\right] \cdot \delta \nu(t_0) &= \left(\frac{\partial T}{\partial \dot q} \circ \pi_{\dot q}\right) \left[\Theta_{t-t_0}\nu(t_0)\right] \cdot \delta \gamma(t)\bigg|_{t_0}^{t_f}\\
&= \left(\frac{\partial T}{\partial \dot q} \circ \pi_{\dot q} \right) \left[\Theta_{t-t_0}\nu(t_0)\right] dq \cdot \delta \nu(t)\bigg|_{t_0}^{t_f}\\
&= \left(\frac{\partial T}{\partial \dot q} \circ \pi_{\dot q} \right) \left[\Theta_{t-t_0}\nu(t_0)\right] dq \cdot {\Theta_{t-t_0}}_* \delta \nu(t_0)\bigg|_{t_0}^{t_f}\\
&= ({\Theta_{t_f-t_0}}^* \theta_L - \theta_L)_{\nu(t_0)} \cdot \delta \nu(t_0),
\intertext{where $\theta_L$ is the one-form $\left(\frac{\partial T}{\partial \dot q}\circ \pi_{\dot q}\right) dq$. Since $dS_i$ is exact,}
d^2S_i &= 0 = {\Theta_{t_f-t_0}}^*d\theta_L - d\theta_L,
\end{align*}
so since $t_0$ and $t_f$ are arbitrary, $\Theta_s^*d\theta_L = d\theta_L$ for arbitrary times $s$, and $\Theta$ preserves the so-called \emph{symplectic form} $d\theta_L$.

\subsection{Discretization \label{s-discrete}}

Discrete mechanics~\cite{aVeselov1988, aSuris1990, aMoser1991, aMarsden2001, bHairer2006} describes a discretization of Hamilton's principle, yielding a numerical integrator that shares many of the structure-preserving properties of the continuous flow $\Theta_s$. Consider a discretization of the trajectory $\gamma: [t_0, t_f] \to \Q$ by a piecewise linear trajectory interpolating $n$ points $\q = \{q_0, q_1, \ldots q_{n-1}\}$, with $q_0 = \gamma(t_0)$ and $q_{n-1} = \gamma(t_f)$, where the discrete velocity $\dot q_{i+1/2}$ on the segment between $q_i$ and $q_{i+1}$ is
\begin{align*}
\dot q_{i+1/2} = \frac{q_{i+1}-q_{i}}{h}, \quad h = \frac{t_f-t_0}{n-1}.
\end{align*}

We seek an analogue of (\ref{eq-EL}) in this discrete setting. To that end, we formulate a discrete Lagrangian
\begin{align}
L_d(q_a, q_b) &= T\left(\frac{q_b-q_a}{h}\right) - V(q_b)\label{eq-DL}
\intertext{and discrete action}
S_d(\q) &= \sum_{i=0}^{n-2} h L_d(q_i, q_{i+1}). \label{eq-daction}
\end{align}

Motivated by (\ref{eq-Hamilton}), we impose a discrete Hamilton's principle:
\begin{align*}
dS_d(\q) \cdot \delta \q = 0
\end{align*}
for all variations $\delta \q = \{\delta q_0, \delta q_1, \ldots, \delta q_{n-1}\}$ that fix $\q$ at its endpoints, \ie, with $\delta q_0 = \delta q_{n-1} = 0.$ For ease of notation, we define versions of the kinetic and potential energy terms in (\ref{eq-DL}) that depend on $(q_a, q_b)$ instead of $(q, \dot q)$:
\begin{align*}
T_d(q_a, q_b) &= T\left(\frac{q_b-q_a}{h}\right) &
T'_d(q_a, q_b) &= \frac{\partial T}{\partial \dot q}\left(\frac{q_b-q_a}{h}\right)\\
V_d(q_a, q_b) &= V(q_b) &
V'_d(q_a, q_b) &= \frac{\partial V}{\partial q}(q_b).
\end{align*}
Then
\begin{align*}
dS_d(\q) \cdot \delta \q &= \sum_{i=0}^{n-2} h \left(D_1 L_d(q_i, q_{i+1}) \cdot \delta q_i + D_2 L_d(q_i, q_{i+1}) \cdot \delta q_{i+1}\right)\\
&= \sum_{i=0}^{n-2} h \left(-\frac{1}{h}T'_d(q_i,q_{i+1}) \cdot \delta q_i + \frac{1}{h}T'_d(q_i, q_{i+1}) \cdot \delta q_{i+1} - \frac{\partial V}{\partial q}(q_{i+1}) \cdot \delta q_{i+1} \right)\\
&= T'_d(q_{n-2},q_{n-1})\cdot \delta q_{n-1} - T'_d(q_0,q_1)\cdot \delta q_0 - h \frac{\partial V}{\partial q}(q_{n-1}) \cdot \delta q_{n-1}\\
&\quad + \sum_{i=1}^{n-2} \left( T'_d(q_{i-1},q_i) - T'_d(q_i,q_{i+1}) - h \frac{\partial V}{\partial q}(q_i)\right)\cdot \delta q_i \\
&= \sum_{i=1}^{n-2} \left( T'_d(q_{i-1}, q_{i}) - T'_d(q_i, q_{i+1}) - h \frac{\partial V}{\partial q}(q_i)\right)\cdot\delta q_i = 0.
\end{align*}
Since $\delta q_i$ is unconstrained for $1 \leq i \leq n-2$, we must have
\begin{align}
\frac{\partial T}{\partial \dot q}(\dot q_{i+1/2}) - \frac{\partial T}{\partial \dot q}(\dot q_{i-1/2}) = -h \frac{\partial V}{\partial q}(q_i), \quad i = 1, \ldots, n-2 \label{eq-DEL},
\end{align}
the \emph{discrete Euler-Langrange equations} of the system.

Unlike in the continuous settings, the discrete Euler-Lagrange equations do not always have a unique solution given initial values $q_0$ and $q_{1}$. We therefore assume in all that follows that $T_d$ and $V_d$ are of a form so that (\ref{eq-DEL}) gives a unique $q_{i+1}$ given $q_i$ and $q_{i-1}$---this assumption always holds, for instance, in the typical case where $T_d$ is quadratic in $\dot q$. Then the discrete Euler-Lagrange equations give a well-defined discrete flow
\begin{align*}
F: (q_{i-1}, q_i) \mapsto (q_i, q_{i+1}),
\end{align*}
which recovers the entire trajectory from initial conditions, in perfect analogy to the continuous setting.

\subsection{Symplecticity of the Discrete Flow}

We now would like a symplectic form preserved by $F$, just as $d\theta_L$ is preserved by $\Theta$. As in the continuous setting, we restrict trajectories $\q$ to those that satisfy the discrete Euler-Lagrange equations, and restrict variations to first variations (and relax the condition that these variations vanish at the endpoints), yielding
\begin{align*}
dS_d(\q) \cdot \delta \q = T'_d(q_{n-2}, q_{n-1})\cdot\delta q_{n-1} - T'_d(q_0,q_1)\cdot \delta q_0 - h \frac{\partial V}{\partial q}(q_{n-1}) \cdot \delta q_{n-1}.
\end{align*}
We denote by $F^k$ the discrete flow $F$ composed with itself $k$ times, or $k$ ``steps'' of $F$. We remark again that all $\q$ satisfying (\ref{eq-DEL}) can be parameterized by initial conditions $\nu_0 = (q_0, q_1)$, and first variations by $\delta \nu_0 = (\delta q_0, \delta q_1)$, so that we can rewrite the discrete action as
\begin{align*}
S_{id}(\nu_0) = \sum_{i=0}^{n-2} h L_d( F^i\nu_0 ).
\end{align*}

Putting together all of the pieces,
\begin{align*}
dS_{id}(\nu_0) \cdot \delta\nu_0 &= dS_d(\q) \cdot \delta \q\\
&= T'_d(q_{n-2},q_{n-1})\cdot \delta q_{n-1} - T'_d(q_0,q_1)\cdot \delta q_0 - h \frac{\partial V}{\partial q}(q_{n-1}) \cdot \delta q_{n-1}\\
&= \left(T'_d(q_a, q_b) - h \frac{\partial V}{\partial q}(q_b)\right) d q_b \cdot (\delta q_{n-2},\delta q_{n-1})\Big\vert_{q_a=q_{n-2},\ q_b=q_{n-1}}\\
&\quad -T'_d(q_a,q_b) dq_a \cdot (\delta q_0, \delta q_1) \Big\vert_{q_a=q_0,\ q_b = q_1}\\
&= \left[T'_d(F^{n-2} \nu_0) - h V'(F^{n-2} \nu_0) \right] d q_b \cdot {F^{n-2}}_* \delta \nu_0 - T'_d(\nu_0)dq_a \cdot \delta \nu_0\\
&= \theta^+_{F^{n-2}\nu_0} \cdot {F^{n-2}}_* \delta \nu_0 + \theta^-_{\nu_0} \cdot \delta\nu_0\\
&= \left({F^{n-2}}^* \theta^+\right)_{\nu_0} \cdot \delta\nu_0 + \theta^-_{\nu_0} \cdot \delta\nu_0.
\end{align*}
for the indicated two-forms $\theta^+$ and $\theta^-$. Since $d(h L_d) = \theta^+ + \theta^-$, $d^2(h L_d) = 0 = d\theta^+ + d\theta^-$. Moreover the intial conditions $\nu_0$ are arbitrary, hence
\begin{align*}
d^2S_{id} = 0 = {F^{n-2}}^* d\theta^+ + d\theta^- = -{F^{n-2}}^*d\theta^- + d\theta^-,
\end{align*}
so
\begin{align*}
d\theta^- = {F^{n-2}}^* d\theta^-.
\end{align*}
Since $n$ is arbitrary, we conclude that the discrete flow $F$ preserves the symplectic form $d\theta^-$. Using backwards error analysis, it can be shown that this geometric property guarantees that integrating with $F$ introduces no energy drift for a number of steps exponential in $h$~\cite{bHairer2006}, a highly desirable property when simulating molecular dynamic or other Hamiltonian systems whose qualitative behavior is substantially affected by errors in energy.

\section{Asynchronous Variational Integrators \label{s-AVI} }

In section \ref{s-discrete} we formulated an action functional
(\ref{eq-daction}) as the integration of a single discrete Lagrangian
over a single time step size $h$. Such a construction is cumbersome
when modeling multiple potentials of varying stiffnesses acting on
different parts of the system: to prevent instability we are forced to
integrate the entire system at the resolution of the stiffest
force. Given a spatial triangulation $\mathcal{T} = \{K\}$ of the
system, asynchronous variational integrators (AVIs), introduced by Lew \etal~\cite{aLew2003}, are a family of numerical integrators, derived from a discrete Hamilton's principle, that support integrating potentials on different triangles at different time steps. In the exposition that follows, we follow the arguments set forth by Lew \etal, but depart at times from the notation used in their work. Although the additional notation and indices introduced herein are initially cumbersome, they will allow for a relatively easy transition to the triangulation-free setting in Section \ref{s-trifree}.

Instead of a global discrete Lagrangian, we instead imbue each triangle $K$ with a local discrete Lagrangian
\begin{align*}
L_d^K(q^K_a, q^K_b) = \int_{t_a}^{t_b} T^K\left[\dot q^K(t)\right]dt - h^K V^K(q^K_b),
\end{align*}
where $T^K$ and $V^K$ are the elemental kinetic and potential energies on triangle $K$, respectively, $h^K = t_b - t_a$ is the elemental time step, and $\dot q^K(t)$, the elemental velocity at time $t$, is left imprecise for the moment. We no longer assume that velocity is constant between times $t_a$ and $t_b$---this would only be true if for every potential on a triangle adjacent to $K$, no multiple of its time step lies between $t_a$ and $t_b$, which is not necessarily the case---so unlike for the discrete Lagrangian (\ref{eq-DL}), here we cannot explicitly integrate the kinetic energy term. For this reason we now write the Lagrangian as an integrated quantity, instead of deferring the integration to inside the action.

Each triangle is only concerned with certain moments in time---namely, integer multiples of $h^K$---and these moments are inconsistent across triangles. We therefore subdivide time in a way compatible with all triangles: for a $\tau$-length interval of time, we define
\begin{align*}
\Xi(\tau) = \bigcup_{K \in \mathcal{T}} \bigcup_{j=0}^{\lfloor \tau/h^K \rfloor} jh^K.
\end{align*}
That is, $\Xi(\tau)$ is the set of all integer multiples less than $\tau$ of all elemental time steps. $\Xi$ can be ordered, and in particular we let $\xi(i)$ be the $(i+1)$-st least element of $\Xi$. If $n$ is the cardinality $\Xi$, we then discretize a trajectory of duration $\tau$ by linearly interpolating intermediate configurations $q_0, q_1, \ldots, q_{n-1}$, where $q_i$ is the configuration of the system at time $\xi(i)$. We discretize velocity as $\dot q_{k+1/2} = \frac{q_{k+1}-q_k}{\xi(k+1)-\xi(k)}$ on the segment of the trajectory between $q_k$ and $q_{k+1}$. We now need to write a global action functional of these trajectories that sums the above elemental Lagrangians, which we do in the natural way:
\begin{align}
S_{\textrm{AVI}}(\q) = \sum_{K \in \mathcal{T}} \sum_{j=0}^{\lfloor \tau/h^K \rfloor} L_d^K\left(q^K_j, q^K_{j+1}\right). \label{eq-actionAVI}
\end{align}

As before, we consider variations $\delta \q = \{ \delta q_0, \ldots, \delta q_{n-1} \}$ with $\delta q_0 = \delta q_{n-1} = 0$, and impose Hamilton's principle,
\begin{align*}
dS_{\textrm{AVI}}(\q) \cdot \delta \q = 0.
\end{align*}
To avoid becoming bogged down in notation, we let $\omega^K(j) = \xi^{-1}(jh^K)$---that is, $\omega$ maps local time indices for $K$ to global indices into $\Xi$---and will write $q_j$ interchangeably for $\pi_K q_j$, the restriction of the (global) configuration $q_j$ to an elemental configuration on $K$. Then
\begin{align*}
S_{\textrm{AVI}}(\q) &= \sum_{K \in \mathcal{T}} \left( \sum_{j=0}^{\lfloor \tau/h^K \rfloor-1} L_d^K\left(q^K_j, q^K_{j+1}\right) + \int_{\lfloor \tau/h^K \rfloor h^K}^{\tau} T^K\left[\dot q^K(t)\right] dt\right)\\
&= \sum_{K \in \mathcal{T}} \left( \sum_{j=0}^{\lfloor \tau/h^K \rfloor-1} L_d^K\left(q_{\omega^K(j)}, q_{\omega^K(j+1)}\right)
+\int_{\lfloor \tau/h^K \rfloor h^K}^{\tau} T^K\left[\dot q^K(t)\right] dt\right)\\
&= \sum_{K \in \mathcal{T}} \left(\sum_{j=0}^{\lfloor \tau/h^K \rfloor-1} \left( \int_{jh^K}^{(j+1)h^K} T^K\left[\dot q^K(t)\right]dt - h^K V^K(q_{\omega^K(j+1)}) \right)+\int_{\lfloor \tau/h^K \rfloor h^K}^{\tau} T^K\left[\dot q^K(t)\right] dt\right)\\
&= \sum_{K \in \mathcal{T}} \left( \sum_{k=0}^{n-2} \left[\xi(k+1)-\xi(k)\right]T^K\left(\frac{q_{k+1} - q_{k}}{\xi(k+1)-\xi(k)}\right) - \sum_{j=0}^{\lfloor \tau/h^K \rfloor-1} h^K V^K(q_{\omega^K(j+1)}) \right).
\end{align*}
Thus, writing
\begin{align*}
{T_d^K}(q_a, q_b, t_a, t_b) &= T^K\left(\frac{q_b-q_a}{t_b-t_a}\right) &
{T_d^K}'(q_a, q_b,t_a,t_b) &= \frac{\partial T^K}{\partial \dot q}\left(\frac{q_b-q_a}{t_b-t_a}\right)\\
V_d^K(q_a, q_b) &= V^K(q_b) &
{V_d^K}'(q_a, q_b) &= \frac{\partial V^K}{\partial q}(q_b),
\end{align*}
we have
\begin{align*}
dS_{\textrm{AVI}}(\q) \cdot \delta \q &= \sum_{K \in \mathcal{T}} \left( \sum_{k=0}^{n-2} {T_d^K}'\left[q_{k},q_{k+1},\xi(k),\xi(k+1)\right]\cdot \left(\delta q_{k+1} - \delta q_{k}\right)\right)  - \sum_{K \in \mathcal{T}} \sum_{j=1}^{\lfloor \tau/h^K \rfloor} h^K \frac{\partial V^K}{\partial q^K}(q_{\omega^K(j)})\cdot \delta q_{\omega^K(j)}\\
&= \sum_{K \in \mathcal{T}} \left( {T_d^K}'\left[q_{n-2},q_{n-1},\xi(n-2),\xi(n-1)\right]\cdot \delta q_{n-1} - {T_d^K}'\left[q_0,q_1,\xi(0),\xi(1)\right]\cdot \delta q_{0}\right)\\
&\quad +\sum_{K \in \mathcal{T}} \sum_{k=1}^{n-2} \left( {T_d^K}'\left[q_{k-1},q_k,\xi(k-1),\xi(k) \right]-{T_d^K}'\left[q_k,q_{k+1},\xi(k), \xi(k+1) \right]\right) \cdot \delta q_{k}\\
&\quad - \sum_{K \in \mathcal{T}} \sum_{j=1}^{\lfloor \tau/h^K \rfloor} h^K \frac{\partial V^K}{\partial q^K}(q_{\omega^K(j)})\cdot \delta q_{\omega^K(j)}\\
&= \sum_{k=1}^{n-2} \sum_{K \in \mathcal{T}} \left({T_d^K}'\left[q_{k-1},q_k,\xi(k-1),\xi(k)\right]-{T_d^K}'\left[q_k,q_{k+1},\xi(k),\xi(k+1)\right]\right) \cdot \delta q_k\\
&\quad - \sum_{k=1}^{n-2} \sum_{h^K \vert \xi(k)} h^K \frac{\partial V^K}{\partial q^K}(q_{k})\cdot \delta q_k,
\end{align*}
where we abuse the notation $h^K \vert m$ to mean, ``all elemental time steps $h^K$ which evenly divide $m$.''
Writing the total kinetic energy of the system $\sum_{K\in \mathcal{T}} T^K$ as $T_{tot}$, for AVIs we recover the discrete Euler-Lagrange equations
\begin{align}
\frac{\partial T_{tot}}{\partial \dot q}(\dot q_{k+1/2}) - \frac{\partial T_{tot}}{\partial \dot q}(\dot q_{k-1/2}) = - \sum_{h^K \vert \xi(k)} h^K \frac{\partial V^K}{\partial q^K}(q_{k}). \label{eq-AVIDEL}
\end{align}
These equations are similar to those we derived for synchronous variational integrators (\ref{eq-DEL}), except that only a subset of potentials $V_d^i$ contribute during each time step. As in the synchronous case, if, as is typical, $T_{tot}$ is quadratic in $\dot q$, the system (\ref{eq-AVIDEL}) gives rise to an explicit numerical integrator that is particularly easy to implement in practice.

\subsection{Multisymplecticity} \label{s-AVImulti}

The right hand side of (\ref{eq-AVIDEL}) depends on $\xi(k)$, and so unlike (\ref{eq-DEL}), the Euler-Lagrange equations for AVIs are time dependent, and do not give rise to a stationary update rule $F(q_{i-1}, q_{i})\mapsto (q_{i}, q_{i+1})$. Instead, we consider the total, time-dependent flow $\hat{F}^k(q_{0}, q_{1}) \mapsto (q_{k}, q_{k+1})$. Once again, we parameterize trajectories satisfying (\ref{eq-AVIDEL}) by $\nu_0=(q_0, q_1)$, and first variations by $\delta \nu_0=(\delta q_0, \delta q_1)$. Restricting ourselves to such trajectories and variations, we rewrite the action (\ref{eq-actionAVI}) as
\begin{align*}
S_{\textrm{iAVI}}=\sum_{K \in \mathcal{T}} \left( \sum_{k=0}^{n-2} \left[\xi(k+1)-\xi(k)\right]T^K_d\left(\hat{F}^k(\nu_0), \xi(k), \xi(k+1)\right) - \sum_{j=0}^{\lfloor \tau/h^K \rfloor-1} h^K V^K_d (\hat{F}^{\omega^K(j+1)-1}(\nu_0)) \right).
\end{align*}

Then
\begin{align*}
dS_{\textrm{iAVI}}(\nu) \cdot \delta \nu &= dS_{\textrm{AVI}}(\q) \cdot \delta \q \\
&= \sum_{K \in \mathcal{T}} \left( {T_d^K}'\left[q_{n-2},q_{n-1},\xi(n-2),\xi(n-1)\right]\cdot \delta q_{n-1} - {T_d^K}'\left[q_0,q_1,\xi(0),\xi(1)\right]\cdot \delta q_{0}\right)\\
&\quad - \sum_{h^K \vert \xi(n-1)} h^K \frac{\partial V^K}{\partial q^K}(q_{n-1})\cdot \delta q_{n-1}\\
&= \sum_{K \in \mathcal{T}} \left( {T_d^K}'\left[\hat{F}^{n-2}(\nu_0), \xi(n-2), \xi(n-1)\right] \cdot \delta q_{n-1} - {T_d^K}'\left[\nu_0, \xi(0), \xi(1)\right]\cdot \delta q_{0}\right)\\
&\quad - \sum_{h^K \vert \xi(n-1)} h^K {V^K_d}'\left[\hat{F}^{n-2}(\nu_0)\right]\cdot \delta q_{n-1}\\
&= \theta^-_{\nu_0} \cdot \delta \nu_0 + \theta^+_{\hat{F}^{n-2}\nu_0} \cdot {\hat{F}^{n-2}}{}_* \delta \nu_0\\
&= (\theta^- + {\hat{F}^{n-2}}{}^* \theta^+)_{\nu_0} \cdot \delta \nu_0
\end{align*}
for one-forms $\theta^-$ and $\theta^+$. Once again we have that
\begin{align}
0 = d^2 S = d\theta^- + \hat{F}^{n-2}{}^*d\theta^+, \label{eq-multi}
\end{align}
but unlike when our action was a sum of Lagrangians, from the \emph{multisymplectic form formula} (\ref{eq-multi}) we have no way of relating $d\theta^-$ to $d\theta^+$, and thus do not recover symplectic structure preservation. Nevertheless, Lew \etal~\cite{aLew2003} conjecture that this multisymplectic structure leads to the good energy behavior observed for AVIs.

\section{Triangulation-Free AVIs} \label{s-trifree}
The above formulation of AVIs assumed a spatial triangulation over which we defined distinct, local Lagrangians. We now present a simple extension that supports potentials with arbitrary, possibly non-disjoint spatial stencil.

Let $\{V^i\}$ be potentials with time steps $h^i$. As in AVIs, for trajectories of duration $\tau$ we define the set of times
\begin{align*}
\Xi(\tau) = \bigcup_{V^i} \bigcup_{j=0}^{\lfloor \tau/h^i\rfloor} jh^i,
\end{align*}
the smallest set of times compatible with the time steps of all of the potentials. Again, let $\Xi$ have cardinality $n$, $\xi(i)$ be the $(i+1)$-th least element of $\Xi$, and $\omega^i(j) = \xi^{-1}(jh^i)$. Then, for $T(\dot q)$ the kinetic energy of the entire configuration, $T_d(q_a, q_b, t_a, t_b) = T\left(\frac{q_b-q_a}{t_b-t_a}\right)$, and $T'_d(q_a, q_b, t_a, t_b) = \frac{\partial T}{\partial \dot q}\left(\frac{q_b-q_a}{t_b-t_a}\right)$, we write the action
\begin{align*}
S_g(\q) = \sum_{j=0}^{n-2} \left[\xi(j+1)-\xi(j)\right]T_d\left[q_j, q_{j+1}, \xi(j), \xi(j+1)\right]- \sum_{V^i} \sum_{j=1}^{\lfloor \tau/h^i\rfloor} h^i V^i(q_{\omega^i(j)}).
\end{align*}
We have made no attempt to define a Lagrangian pairing the kinetic and potential energy terms; we will see that an action defined this way still leads to a multisymplectic numeric integrator.

To that end we impose $dS_g(\q) \cdot \delta \q = 0$ for variations with $\delta q_0 = \delta q_{n-1} = 0$. Then we rewrite $S_g$ as
\begin{align*}
S_g(\q) = \sum_{j=0}^{n-2} \left[\xi(j+1)-\xi(j)\right]T_d\left[q_j, q_{j+1}, \xi(j), \xi(j+1)\right]- \sum_{j=1}^{n-1} \sum_{h^i \vert \xi(j)} h^i V^i(q_j)
\end{align*}
so that
\begin{align*}
dS_g(\q) \cdot \delta \q &= \sum_{j=0}^{n-2} T'_d\left[q_j, q_{j+1}, \xi(j), \xi(j+1)\right]\cdot \left(\delta q_{j+1} - \delta q_j\right) - \sum_{j=1}^{n-1} \sum_{h^i \vert \xi(j)} h^i \frac{\partial V_i}{\partial q}(q_j)\cdot \delta q_j\\
&= T'_d\left[q_{n-2}, q_{n-1}, \xi(n-2), \xi(n-1)\right]\cdot \delta q_{n-1} - T'_d\left[q_0, q_1, \xi(0), \xi(1)\right]\cdot \delta q_0\\
&\quad - \sum_{h^i \vert \xi(n-1)}h^i \frac{\partial V^i}{\partial q}(q_{n-1}) \cdot \delta q_{n-1}\\
&\quad + \sum_{j=1}^{n-2} \left( T'_d\left[q_{j-1}, q_{j}, \xi(j-1), \xi(j)\right] - T'_d\left[q_j, q_{j+1}, \xi(j), \xi(j+1)\right]-\sum_{h^i \vert \xi(j)} h^i \frac{\partial V^i}{\partial q}(q_j)\right)\cdot \delta q_j\\
&= \sum_{j=1}^{n-2} \left( T'_d\left[q_{j-1}, q_{j}, \xi(j-1), \xi(j)\right] - T'_d\left[q_j, q_{j+1}, \xi(j), \xi(j+1)\right]-\sum_{h^i \vert \xi(j)} h^i \frac{\partial V^i}{\partial q}(q_j)\right)\cdot \delta q_j.
\end{align*}
The Euler-Lagrange equations are then
\begin{align}
\frac{\partial T}{\partial \dot q}(\dot q_{k+1/2}) - \frac{\partial T}{\partial \dot q}(\dot q_{k-1/2}) = - \sum_{h^i \vert \xi(k)} h^i \frac{\partial V^i}{\partial q^i}(q_{k}), \label{eq-gAVIDEL}
\end{align}
exactly the same as the Euler-Lagrange equations (\ref{eq-AVIDEL}) for
ordinary AVIs. Triangulation-free AVIs can thus be integrated in
exactly the same manner as ordinary AVIs, for instance, by using the
algorithm presented by Lew et al.~\cite{aLew2003}.

\subsection{Multisymplecticity}
To show that triangulation-free AVIs still satisfy the multisymplectic form formula (\ref{eq-multi}), we follow the derivation for multisymplecticity of ordinary AVIs. Replacing $\sum_{K \in \mathcal{T}} T_d^K$ with $T_d$ in Section \ref{s-AVImulti}, an identical calculation shows triangulation-free AVIs satisfy (\ref{eq-multi}).

\section{Sphere-plate Impact}\label{sec:results}

Our triangulation-free multisymplectic formulation supports interaction potentials of the kind needed for contact response, and we expect any such method to exhibit the good energy behavior associated with multisymplectic integrators. As a numerical experiment of this behavior, we simulated the impact of a spherical shell with a thin plate, as described in Cirak and West's article on Decomposition Contact Response (DCR)~\cite{aCirak2005}, using the Asynchronous Contact Mechanics (ACM) framework~\cite{aHarmon2009} built on triangulation-free AVIs. A sphere of radius 0.125 approaches a plate of radius 0.35 with relative velocity 100. Both the sphere and the plate have thickness 0.0035.
The time steps of our material forces (stretching and bending) are $10^{-7}$ (the same as those chosen by Cirak and West.)

\begin{figure}[t]
\begin{center}
\includegraphics[width=4in]{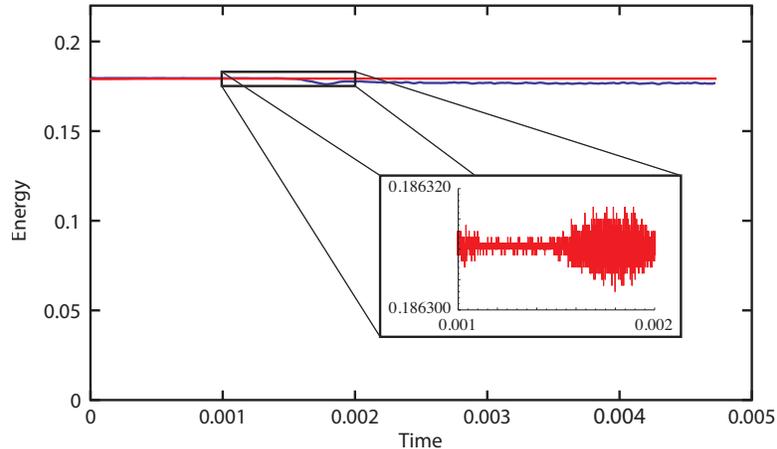}
\caption{Total energy over time of a thin sphere colliding against a
  thin plate, simulated using asynchronous contact mechanics~\cite{aHarmon2009} (red) compared to data provided for decomposition
  contact response~\cite{aCirak2005} (dark blue).}
\label{fig-plate}
\end{center}
\end{figure}

Figure \ref{fig-plate} compares energy over time when this simulation is run using both our ACM and DCR.
Using ACM there is no noticeable long-term drift.
Closely examining the energy data produced by ACM reveals the
high-frequency, low-amplitude, qualitatively-negligible oscillations
characteristic of symplectic integrators.

\paragraph{Acknowledgements}
We thank Fehmi Cirak for providing the meshes used in the comparison
to DCR.  This work was supported in part by the NSF (MSPA Award
No. IIS-05-28402, CSR Award No. CNS-06-14770, CAREER Award
No. CCF-06-43268). The Columbia authors are supported in part by
generous gifts from Adobe, ATI, Autodesk, mental images, NVIDIA, the
Walt Disney Company, and Weta Digital.

\bibliography{refs}

\begin{thebibliography}{10}

\bibitem{aCirak2005}
F.~Cirak and M.~West.
\newblock Decomposition contact response ({DCR}) for explicit finite element
  dynamics.
\newblock {\em International Journal for Numerical Methods in Engineering},
  64(8):1078--1110, 2005.

\bibitem{bHairer2006}
E.~Hairer, C.~Lubich, and G.~Wanner.
\newblock {\em Geometric Numerical Integration: Structure-Preserving Algorithms
  for Ordinary Differential Equations}, volume~31 of {\em Springer Series in
  Computational Mathematics}.
\newblock Springer-Verlag, second edition, 2006.

\bibitem{aHarmon2009}
D.~Harmon, E.~Vouga, B.~Smith, R.~Tamstorf, and E.~Grinspun.
\newblock Asynchronous contact mechanics.
\newblock In {\em SIGGRAPH '09: ACM SIGGRAPH 2009 papers}, pages 1--12, New
  York, NY, USA, 2009. ACM.

\bibitem{bLanczos1986}
C.~Lanczos.
\newblock {\em The Variational Principles of Mechanics}.
\newblock Dover Publications, fourth edition, 1986.

\bibitem{tLew2003}
A.~Lew.
\newblock {\em Variational Time Integrators in Computational Solid Mechanics}.
\newblock PhD thesis, California Institute of Technology, 2003.

\bibitem{aLew2003}
A.~Lew, J.E. Marsden, M.~Ortiz, and M.~West.
\newblock Asynchronous variational integrators.
\newblock {\em Arch. Rational Mech. Anal.}, 167(2):85--146, 2003.

\bibitem{cMacKay1992}
R.~MacKay.
\newblock Some aspects of the dynamics of hamiltonian systems.
\newblock In {\em The Dynamics of Numerics and the Numerics of Dynamics}. 1992.

\bibitem{aMarsden2001}
J.~Marsden and M.~West.
\newblock Discrete mechanics and variational integrators.
\newblock {\em Acta Numerica}, 10:357--514, 2001.

\bibitem{aMoser1991}
J.~Moser and A.~Veselov.
\newblock Discrete versions of some classical integrable systems and
  factorization of matrix polynomials.
\newblock {\em Communications in Mathematical Physics}, 139(2):217--243, 1991.

\bibitem{aSuris1990}
Y.~Suris.
\newblock Hamiltonian methods of runge-kutta type and their variational
  interpretation.
\newblock {\em Math. Modelling}, 2(4):78--87, 1990.

\bibitem{aVeselov1988}
A.~Veselov.
\newblock Integrable discrete-time systems and difference operators.
\newblock {\em Functional Analysis and Its Applications}, 22(2):83--93, 1988.

\bibitem{pWest2004}
M.~West.
\newblock {\em Variational Integrators}.
\newblock PhD thesis, California Institute of Technology, 2004.

\end{thebibliography}
\bibliographystyle{plain}
\end{document}